%% file: packings.tex
\documentclass[twoside,openright,11pt,a4paper,epsf]{article}
\usepackage{latexsym}
\usepackage[T1]{fontenc}
\usepackage{amsmath}
\usepackage{amsthm}
\usepackage{amsfonts} 
\usepackage{amssymb}
\usepackage{vmargin}
\usepackage{epsfig}
\usepackage[all]{xy}

\newcommand{\color}[6]{}

\newcommand{\R}{\mathbb{R}}
\newcommand{\N}{\mathbb{N}}
\newcommand{\Z}{\mathbb{Z}}
\newcommand{\C}{\mathbb{C}}
\newcommand{\D}{\mathbb{D}}

\renewcommand{\P}{\mathbb{P}}

\newcommand{\ca}{\mathbb{\mathcal A}}

\newcommand{\cc}{\mathbb{\mathcal C}}

\renewcommand{\sc}{\mathcal S}
\newcommand{\ce}{\mathcal E}
\newcommand{\cu}{\mathcal U}

\newcommand{\fonction}[5]
{$$ 
\begin{array}{rcccl}
 #1 & : & #2 & \longrightarrow &#3 \\
    &   & #4 & \longmapsto &#5 
\end{array}
$$}
\renewcommand{\span}{\text{Span}}

\newcommand{\im}{\text{Im}\,}

\newcommand{\priv}{\backslash}
\newcommand{\lra}{\longrightarrow}
\newcommand{\hra}{\hookrightarrow}
\newcommand{\ssi}{\Longleftrightarrow}
\newcommand{\om}{\omega}
\newcommand{\eps}{\varepsilon}
\renewcommand{\phi}{\varphi}

\renewcommand{\L}{\mathcal{L}}

\newcommand{\wdt}[1]{\widetilde{#1}}
\newcommand{\Span}{\text{Span}}
\newcommand{\cqfd}{\hfill $\square$ \vspace{0.1cm}\\ }

\newcommand{\ds}{\displaystyle}

\newtheorem{definition}{Definition}[section]
\newtheorem{lemma}[definition]{Lemma}
\newtheorem{thm}{Theorem}

\newtheorem{cor}[definition]{Corollary}
\newtheorem{prop}[definition]{Proposition}
\newtheorem{rk}[definition]{Remark}

\title{\vspace{-2cm} Maximal symplectic packings in $\P^2$.}
\author{Emmanuel Opshtein
\footnote{This research was supported by THE ISRAEL SCIENCE FOUNDATION 
(grant No. 1227/06).}
}\date{}
\begin{document}
\maketitle\vspace{-.7cm}
\begin{abstract}
In this paper we describe the intersection between the balls of 
maximal symplectic packings of $\P^2$. This analysis shows the existence of 
singular points for maximal packings of $\P^2$ by more than three equal balls.
It also yields a construction of a class of very regular examples of 
maximal packings by five balls. 
\end{abstract}
\section{Introduction.}
The symplectic packing problem is the question of identifying the
conditions on the radii of $k$ balls for being able to pack them
symplectically in a given manifold. It was first considered by Gromov
as a problem whose answer singularizes symplectic geometry from the
volume-preserving one  : the restrictions are sometimes stronger than
the volume obstruction alone \cite{gromov}. For instance, the non-squeezing 
theorem asserts that no ball of radius bigger than one can be packed
(or embedded in this single ball situation) in the infinite volume
cylinder $B^2(1)\times \R^{2n}$. 
In the special case of $\P^2$, this problem has been given a complete
answer. Apart from the volume obstruction, symplectic packings by less
than eight balls are submitted to finitely many purely symplectic
obstructions discovered by Gromov \cite{gromov} and McDuff-Polterovich
\cite{mcpo}. Biran also proved that the symplectic obstructions
disappear for more than nine balls \cite{biran2}. 

This paper is aimed at describing what the packings of $\P^2$ by balls of
maximal radii look like. We are particularly interested in understanding 
their intersection properties. Let us first define our objects. 
\begin{definition}
A maximal symplectic packing of $M$ by $k$ balls is a symplectic embedding
$\phi:(B(r_1)\coprod\dots\coprod B(r_k),\om_{\text{st}})\hra (M,\om)$
where the radii are such that there exists no symplectic packing of
$M$ by balls of radii $(r_1,\dots,r_i+\eps,\dots,r_k)$. 
It will be said {\it smooth} if each $\phi_i:=\phi_{|B(r_i)}$ extends
to a smooth embedding of the closed balls in $M$.
It will be said {\it regular} if these maps have only a finite number of
singular points on the boundary : each $\phi_i$ extends to a
topological embedding of the closed ball which is a 
smooth embedding of $\overline{B(r_i)}\priv\{p_i^1,\dots,p_i^{n_i}\}$.  
\end{definition}
\noindent The space of regular maximal symplectic packings by a fixed number of balls is
naturally endowed with the Hausdorff topology for compact sets (for
instance). Throughout this paper, genericity is meant with respect to
this topology, and should be understood in a strong sense : a
property is generic if it is true for an open dense set. 
Our first theorem deals with smooth maximal symplectic packings of $\P^2$.
\begin{thm}\label{smooth} Below are the generic pattern of
  intersection between the balls of smooth maximal symplectic packings of $\P^2$. \\
a. Generically, the closed balls of a smooth maximal symplectic
packing of $\P^2$ by two balls intersect precisely along one common
Hopf circle of their boundary.\\
b. Generically, any two closed balls of a smooth maximal
symplectic packing of $\P^2$ by three equal balls intersect precisely
along one common Hopf circle of their boundary.\\
c. Generically, the two smallest balls of a smooth
maximal symplectic packing of 
$\P^2$ by three non-equal balls do not intersect, while the intersection of the
biggest ball with any of the others is exactly one common Hopf
circle of their boundaries. \\   
d. There exist no smooth maximal symplectic packing of $\P^2$ by more 
than three equal balls.
\end{thm}
\noindent The Hopf circles of a ball are the characteristic leafs of its boundary 
(see section \ref{partsul}).  Theorem \ref{smooth} is  
better understood in the light of Karshon's examples \cite{mcpo}: it says that 
generic smooth maximal symplectic packing look very much like those she produced  
(see section \ref{examples}). The idea consists in translating 
the maximality property to the existence of characteristic circles in the
intersections of the boundary spheres.
These characteristics give rise to (maybe singular) symplectic spheres, whose
intersection properties lead to the desired uniqueness.  The approach
is based on a strong connection  observed by Paternain-Polterovich-Siburg 
\cite{popasi} or Laudenbach-Sikorav \cite{lasi} between symplectic non-removable
intersection and closed characteristics.

The importance of the characteristic foliation in the present study is
precisely the reason for our definition of smooth or regular maximal
symplectic packing to be so restrictive. Before going further, a
discussion about  the existence of the objects under consideration is 
needed.  A non-contructive argument due to McDuff shows that there
always exist symplectic packings by {\it open} balls of maximal radii,
with no guaranteed boundary regularity. In the other hand,
several explicit examples are available.  We already mentionned Karshon's 
smooth packings by two or three balls. Generalizing her construction, Traynor
\cite{traynor} and Schlenk \cite{schlenk}  
produced examples of maximal symplectic packings of $\P^2$ by five and
six balls, which unfortunately fail from far to be
regular. It is not completely surprising in view of
theorem \ref{smooth}.d.  As far as we are concerned, the achieved boundary regularity
is  nevertheless as difficult to handle as McDuff's abstract maximal 
packings : no convenient notion of characteristic foliation on
the boundary of the balls can be defined. The second result of this
paper concerns precisely the relevance of our definition of regularity.
Allowing only finitely many singularities enables us to produce 
interesting maximal packings of $\P^2$, at least by five balls. 
\begin{thm}\label{ex5}
There exist regular maximal symplectic packings of $\P^2$ by five equall balls.  
\end{thm}
\noindent The construction relies on a decomposition theorem of K\"ahler manifolds 
due to Biran \cite{biran3}.

Our third result generalizes theorem 1 to the regular setting.
The intersections between the balls of a regular maximal
symplectic packing is a union of {\it Hopf
circles} of the boundary spheres which provides a ``grid'' of a
topological $2$-fold, symplectic away from its singularities. 
\begin{thm}\label{regular}
Given a regular maximal symplectic packing of a
symplectic manifold, there exists at least one ``supporting surface'' : 
a closed topological surface covered by the balls of the packing and 
whose intersection with any ball is a union of smooth symplectic discs bounded 
by Hopf circles. Generically, the balls of the packing intersect exactly along the
Hopf circles contained in the supporting surfaces. 
\end{thm}
The existence result above can be sharpened when the singularities are simple 
enough  (see definition \ref{verysimple}). For five balls for instance, the 
intersection pattern must be the same as in the constructed examples (theorem \ref{ex5}).
\begin{thm}\label{theend}
Regular maximal symplectic packings of $\P^2$ by five equal balls which have simple type
have exactly one supporting surface, of symplectic area $2\pi$, intersecting each ball
through exactly one Hopf disc. Generically these maximal packings thus intersect along 
exactly one Hopf circle of each of the boundary spheres. 
\end{thm}
The paper is organized as follows. In section \ref{examples}, we
discuss previously known examples of maximal symplectic packings
and we construct new ones (theorem \ref{ex5}). We hope this section to
shed light on the statements of theorem \ref{smooth}, \ref{regular} 
and \ref{theend} by providing relevant illustrations. In the
third section, we explain
the link between non-removable intersections and characteristic
foliations in the setting of smooth balls, and prove theorem
\ref{smooth}. The purpose of section \ref{partregularization} is to
adapt to non-smooth objects the tools we use in the preceding
section. We prove theorem \ref{regular} in section 
\ref{partreg} and conclude by a technical paragraph   
aimed at smoothening the supporting surfaces in view of proving theorem
\ref{theend}.\vspace{.2cm}

\noindent{\bf Aknowledgements} I would like to express gratitude to P. Biran,
L. Polterovich and J.Y. Welschinger for showing interest to this work, for 
usefull and motivating discussions.  
\section{Examples of maximal symplectic packings.}\label{examples}
\begin{figure}[h!]
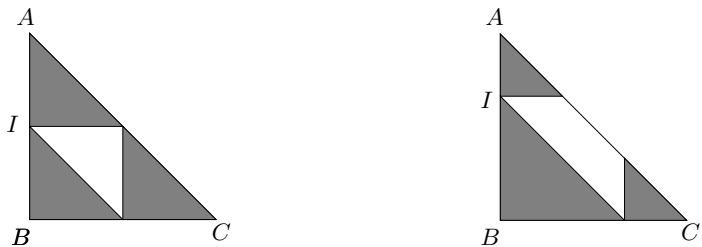

\begin{center} 
\input karshon.pstex_t
\caption{Maximal symplectic packings of $\P^2$ by $3$ balls.}
\label{figkar}
\end{center}\vspace{-.7cm}
\end{figure}

The aim of this part is to provide the reader with examples of smooth
or regular symplectic packings. We describe them in the light of
the results stated in the above introduction. \vspace{.2cm}\\
\textbf{Karshon's construction.}
Examples of smooth maximal symplectic packing of $\P^2$ were shown by
Karshon \cite{mcpo} or Traynor \cite{traynor}. They can be described
in the following way. The momentum map - or the action-angle
coordinates - presents $\P^2$ as a singular
bundle over a closed  triangle with $2$-dimensional tori as generic
fibers. The balls forming the packing are fibered by these tori and
project by the momentum map to close triangles (see figure \ref{figkar}).

\noindent Note that the intersection between the closed balls in these
examples are precisely one common circle of Hopf fibration of the
boundary spheres (point $I$). Moreover, the fiber over the line $AB$
is a $2$-dimensional sphere covered by the union of the two balls
$B_1$, $B_2$, and whose intersection with each ball is a disc bounded
by a Hopf circle. It is thus precisely one of the ``supporting
surfaces'' of the packing. Theorem \ref{smooth} says that every maximal
packings  of $\P^2$ by two or three balls always have the same
intersection patterns as those examples.

From this example, it is also easy to see that the assumption concerning
the equality of the balls in theorem \ref{smooth}.d cannot be
dropped. Actually, if $B_1$, $B$ is a smooth packing of $\P^2$ by 
balls of radii $r\ll 1$ and $r_1=\sqrt{1-r^2}$, there exist $k-1$ 
transformations of $\P U(2)$ leaving $B_1$ invariant and taking $B$ 
to  disjoint symplectic balls $B_2$,\dots, $B_k$. Together with $B_1$, 
they provide a maximal packing of $\P^2$ by $k$ balls.\vspace{.2cm}\\
\textbf{Biran's decomposition theorem (see \cite{biran3}).}  
The examples of regular packings to come are based on a decomposition result
due to Biran which we describe now briefly. Given an integral
K\"ahler manifold $(M,\om,i)$ ({\it i.e.} with $\om\in H^2(M,\Z)$),
there always exist a complex hypersurface $X$ which is Poincare dual
to $k\om$ for some big enough integer multiplier $k$. Biran
showed that there exists a skeleton $\Delta_X\subset M$ (of empty
interior) associated to this hypersurface, whose 
complement in $M$ is a {\it standard symplectic bundle over $X$}. 
It turns out that some maximal packings appear very clearly
in these special coordinates. In order to explain this point, 
let us discuss briefly the structure of these standard bundles. 
Their symplectic type is that of the unit disc bundle associated to a 
Hermitian line bundle of first Chern class $c_1=[k\tau]$, where
$\tau=\om_{|X}$. The symplectic stucture itself is explicitly given by
$$
\om:=\pi^* \tau + d(r^2\alpha)
$$
where $r$ is the radial coordinates in the fiber and $\alpha$ is the
$1$-form whose restriction to the fibers is 
$\alpha_{|\pi^{-1}(x)}=1/k\,d\theta$ on $E\priv\{r=0\}$ and $d\alpha=-\pi^* \tau$.
Notice that although $\theta$ is not defined globally on $E$ because
of non-vanishing Chern class, the differential form ``$d\theta$'' is
perfectly defined out of the zero section.

The connection between standard symplectic bundles and balls result from
the following simple observation. The restriction of these bundles to a
symplectic ellipso\"id of $X$ is symplectomorphically equivalent to an
 ellipso\"id. In the following lemma, 
$\ce_a=\ce_{(a_1,\dots,a_n)}$ denotes the standard ellipso\"id of
$\C^n$ :
$$
\ce_{a}:=\left\{(z_1,\dots,z_n)\in \C^n\; | \;
\frac{|z_1|^2}{\pi a_1}+\dots+\frac{|z_n|^2}{\pi a_n}<1\right\}.
$$
\begin{lemma}\label{cse} Consider the trivial disc bundle
$\pi:\ce_a\times \D\lra\ce_a$ over an ellipso\"id of $\C^n$, equipped
with a symplectic structure defined by $\om:=\pi^*\om_{\text{st}}+d(r^2\alpha)$,
$\alpha_{|\{x\}\times \D}=1/k d\theta$ and $d\alpha=-\pi^*\om_{\text{st}}$,
where $(r,\theta)$ are the polar coordinates on $\D$. Then there
exists a smooth function $h:\ce_{a}\lra
  \R$ such that the map \fonction{\Phi}{(\ce_a\times 
  \D,\om)}{(\ce_{a,\frac{1}{k}},\om_{st})}{(z,w)}{(\sqrt{(1-|w|^2)}z,\frac{1}{\sqrt k}
   e^{i\,h(z)}w)} is a symplectomorphism.
\end{lemma}
\noindent\underline{Proof :} Consider the  coordinates
$(z,w)$=$(r_1,\theta_1,\dots,r_n,\theta_n,r,\theta)$ on $\ce_a\times
\D$. The symplectic form $\om$ is given in these 
coordinates by 
$$
\om=\sum_{i=1}^n dr_i^2\wedge d\theta_i+d(r^2\alpha).
$$
Taking into account the identities  $\alpha_{|\{x\}\times \D}=1/k\,
d\theta$ and $d\alpha=-\sum dr_i^2\wedge d\theta_i$, we get :
$$
\begin{array}{lcl}
\om &= & (1-r^2)\sum_{i=1}^n dr_i^2\wedge d\theta_i+dr^2\wedge \alpha\\
& = & \sum_{i=1}^n d[(1-r^2)r_i^2]\wedge d\theta_i+dr^2\wedge
    [\alpha+\sum_{i=1}^n r_i^2d\theta_i]\\
& = &\ds  \sum_{i=1}^n d[(1-r^2)r_i^2]\wedge d\theta_i+\frac{1}{k}dr^2\wedge[d\theta+\beta],
\end{array}
$$
where $\beta:=k[\alpha-1/k\,d\theta+\sum_{i=1}^n r_i^2d\theta_i]$. 
The form $\beta$ is defined on $\ce_a\times \left(\D\priv
\{0\}\right)$ and is closed. Moreover, its action on
the fundamental group of $\ce_a\times \D\priv\{0\}$ is trivial because
$\beta_{|\{x\}\times \D}=0$. Hence there is 
a smooth function $h:\ce_a\times \D\priv \{0\}\lra \R$ such that
$\beta =dh$. Notice now that $h$ does not depend on $w$
because $\beta$ vanishes on the vertical discs. This function thus
extends to $\ce_a\times \D$ and depends only on $z$.  We finally get :
\begin{eqnarray*}
\om &= & \sum_{i=1}^nd[(1-r^2)r_i^2]\wedge d\theta_i + \frac{1}{k}dr^2\wedge
d(\theta+h)\\
 & =& \Phi^*\left[\sum_{i=1}^n dr_i^2\wedge d\theta_i+dr^2\wedge d\theta\right],
\end{eqnarray*}
where $\Phi$ is the announced map. It clearly sends $\ce_a\times \D$ to
the ellipso\"id $\ce_{a,\frac{1}{k}}$.\cqfd
Since a complex hypersurface of a K\"ahler manifold is K\"ahler, an
obvious iteration leads to the following corollary. It seems to hold true 
also in general compact symplectic manifolds due to Donaldson's results on the 
existence of symplectic hypersurfaces \cite{donaldson}. 
\begin{cor} Every K\"ahler manifold has full packing by one ellipso\"id.
\end{cor}
\noindent\textbf{Maximal full packing of $\P^n$ by $k^n$ balls.} We make a
digression at this point to explain how to construct a full symplectic
packing of $\P^2$ by four {\it open} equal balls of radius $1/\sqrt 2$.
The generalization to $k^n$ balls of maximal 
radius $1/\sqrt[n] k$ in $\P^n$ is straightforward. As far as I know,
although their existence is well known from McDuff-Polterovich's work,
no example of such packings was available up to now.

\begin{figure}[h!]
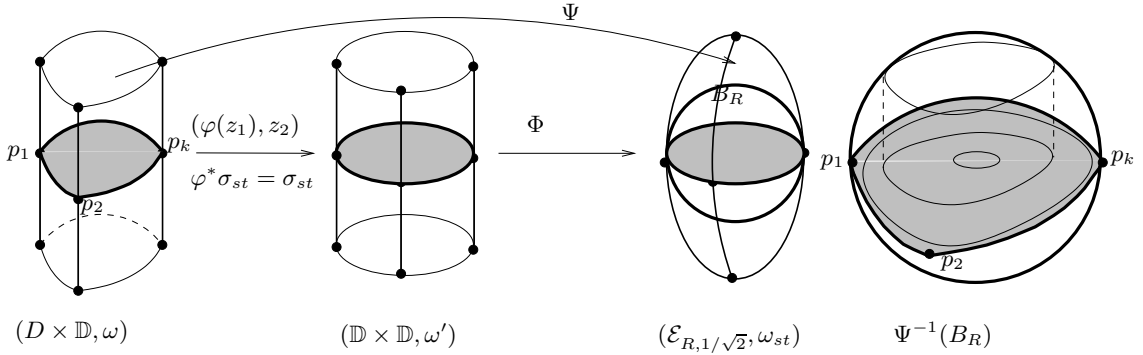

\begin{center} 
\input pack5.pstex_t
\caption{Regular packing of a standard disc bundle over a disc by a ball.}
\label{pack5}
\end{center}
\end{figure}

Consider the quadric $Q:=\{z_0^2+z_1^2+z_2^2=0\}$ in $\P^2$ with
homogeneous coordinates $[z_0:z_1:z_2]$. Then $\pi:\P^2\priv \R\P^2\lra Q$ 
is a standard disc bundle with fibers of area $\pi/2$. Divide $Q$ in four 
open discs $D_i$ of area
$\pi/2$, $i=1,\dots,4$. Then the sets $B_i:=\pi^{-1}(D_i)\subset
(\P^2,\om_{FS})$ are obviously disjoint. Moreover, the previous lemma
shows that they are  symplectic balls of radius $1/\sqrt 2$. Notice
also that the boundary singularities of these balls (inavoidable by
theorem \ref{smooth}.d) are easily descriptible in terms of the
singularities of the discs $D_i$ thanks to the explicit formula for
the symplectomorphism $\Phi$ of lemma \ref{cse}.  \vspace{.2cm}\\
\textbf{Regular maximal packing of $\P^2$ by five balls (see figure
  \ref{pack5}).} As above consider the quadric $Q:=\{z_0^2+z_1^2+z_2^2=0\}$, 
and the projection $\pi:\P^2\priv \R\P^2\lra Q$ which gives $\P^2\priv \R\P^2$
the structure of a
standard disc bundle with fibers of area $\pi/2$. Cover $Q$ by five
closed discs of area $2\pi/5$ with finite number of singularities on
their boundaries. We claim that we can find the desired balls of the
packing inside the saturated sets upon these discs.
To see this, we construct a regular symplectic embedding of a ball of
radius $R$ inside $\pi^{-1}(D)$, where $D$ is any closed disc in $Q$
of area $\pi R^2<\pi/2$, with a finite number of singularities
$p_1,\dots,p_k\in \partial D$. First identify
$\pi^{-1}(D)$ with $D\times \D$. Then $(\pi^{-1}(D),\om_{FS})$ is
symplectomorphic to a standard ellipso\"id $\ce:=\ce_{R,1/\sqrt
  2}$ {\it via} a map $\Psi$ which is the composition of a
fibered map sending $D\times \D$ to the standard bidisc $\D_R\times \D$
with the map $\Phi$ of lemma \ref{cse}. The boundary regularity of the
symplectomorphism $\Psi:\pi^{-1}(D)\lra \ce$ can be
easily described. First it extends to a homeomorphism between
$\overline{D}\times \D$ ({\it i.e.} $\overline{\pi^{-1}(D)}$ minus the
section at infinity) and $\overline{\ce}\priv C_\infty$
where $C_\infty:=\{|z_2|^2=1/\sqrt 2, \; z_1=0\}$. Moreover, this
extension is a local diffeomorphism except at the singular
points $\cup \{p_i\}\times \D$ of $\pi^{-1}(D)$. Now the ellipso\"id
$\overline{\ce}$ contains an euclidean closed ball $B$
of radius $R$ whose boundary intersects $\partial \ce$
only along the ``zero section'' $\{z_2=0,\; |z_1|=R\}$. The map
$\Psi^{-1}:B\lra \pi^{-1}(D)$ is therefore a regular symplectic
embedding of a ball of radius $R$ in $\pi^{-1}(D)$.\cqfd

\indent It may be worth noticing that this construction confirms the 
intuition that the space of maximal symplectic packings of $\P^2$ by four 
or five balls is not connected (in contrast with the non-maximal situation,
see \cite{biran1,mcduff}). It is less clear for two or three balls.   
\section{Non-removable intersections in smooth maximal
  packings.}
The aim of this section is to introduce the main tool of the paper -
namely the link between non-removable intersections in symplectic
geometry and characteristic foliations - following \cite{popasi,lasi}. 
We also prove theorem \ref{smooth}. Dealing  only with smooth objects 
first permits to avoid  the technical difficulties arising in
the context of regular packings. The ideas should be more transparent
to the reader.  \\
\indent Let us explain first how the
genericity is achieved in theorem \ref{smooth}. 
The Hamiltonian transformations of $\P^2$ act trivally on the space of
maximal symplectic packings by moving all balls together. More
generally, one can also move on this space by considering independant
Hamiltonian transformation of each ball of the packing, as long as
they preserve the packing property (the interior of the
balls do not intersect). Precisely, we define a {\it Hamiltonian
  perturbation of a maximal packing} $\{\phi_i\}_{i=1,\dots,k}$ as a family
$\{\phi_i^t\}:=\{\Phi^t_{X_{H_i}}\circ\phi_i\}$, where the
time-dependent Hamiltonian functions $H_i$ on $M$ are such that the 
open balls $\phi_i^t(B(r_i))$ remain disjoint sets for all $t\in [0,\eps]$. Note that, 
unlike global Hamiltonian transformations, these Hamiltonian perturbations 
of packings enable to break intersections. The properties stated 
in theorem \ref{smooth} are generic precisely because they are always true 
after a possible Hamiltonian perturbation of the packing.
\subsection{Digging in a ball by a Hamiltonian : Sullivan's lemma.}\label{partsul}
As explained in  \cite{popasi} and \cite{lasi}, the characteristic
foliation plays a central role in the phenomena of ``symplectic
non-removable intersection''. Loosely speaking, the reason for a
compact set inside the smooth boundary of an open set $U$ not to be
displacable inside $\overline{U}$ by a Hamiltonian vector field is
that it contains a closed invariant set of the characteristic
foliation of $\partial U$. 

\begin{definition} The characteristic
  distribution $\L$ of a hypersurface $S$ in a symplectic manifold
  $(M,\om)$ is the kernel of the restriction of $\om$ to $S$
  ({\it i.e.} $\forall x\in S, \; \L_x:=\ker {\om_x}_{|T_xS}$). The
  characteristic foliation of $S$ is the integral foliation 
associated to this one-dimensional characteristic distribution. 
\end{definition}
The characteristic foliation is obviously preserved by any map 
$\phi:S\subset (M,\om)\lra S'\subset(M',\om')$ with 
$\phi^*\om'_{|TS'}=\om_{|TS}$.\vspace{.1cm}\\
\textbf{Example :} On the euclidean sphere $S(r)\subset \C^n$,
the characteristic distribution is given by $\L_x=\text{Span}_\R(i\vec
N(x))$ where $\vec N(x):=x/\Vert x\Vert$. The characteristic foliation
is thus the classical foliation of $S(r)$ by Hopf
circles.
We call {\it smooth (or regular) symplectic closed ball} of $(M,\om)$ any symplectic
 smooth (or regular) embedding of a euclidean closed ball of $\C^n$ in $M$.
If $B$ is a symplectic closed ball in $M$
corresponding to a regular embedding $\phi:\overline{B(r)}\lra M$, the
characteristic leaves of $\partial B=\phi(S(r))$ are the images by
$\phi$ of the Hopf circles of $S(r)$ and will be called the Hopf
circles of $B$. Finally, the Hopf discs of a symplectic ball
$B$ are the images of the intersections of $B(r)$ with the complex lines
of $\C^n$.\vspace{.1cm}

The vector field $i\vec N(x)$ defines and orients the characteristic foliation
of the euclidean sphere $S$ in $\C^n$. The Hamiltonian vector
field $\vec X_H$ associated to a smooth function $H:\C^n\lra \R$ points inside
$B$ at a point $x$ of $S$ if and only if $i\vec N\cdot H (x)<0$. Actually,
$$
\vec X_H\cdot \vec N (x)<0\ssi \om(\vec X_H(x),i\vec N(x))<0\ssi
dH_x(i\vec N(x))<0.
$$
In other terms, the Hamiltonian vector field $\vec X_H$ ``digs'' in
$B$ at a point $x$ of $S$ if {\it $H$ is decreasing along the 
  characteristic foliation at this point}. The following lemma explains along
which compact sets of $S$ one can dig in $B$ in a Hamiltonian way. It
is a particular case of a more general result due to Sullivan
\cite{sullivan} (see also \cite{lasi}) :
\begin{lemma} \label{sul} Let $K$ be a compact set of $S^{2n-1}\subset \C^n$ and
  $\cc$ the set of Hopf circles belonging to $K$. There exists a
  smooth function $H:S^{2n-1}\lra \R$ such that $i\vec N\cdot H<0$ on
  $K\priv \cc$ and $dH=0$ on $\cc$. 
\end{lemma} 
\noindent \underline{Proof :} We first reduce Lemma \ref{sul} to
finding convenient functions on solid tori by an argument of partition
of unity. Denote $\pi:S^{2n-1}\lra \P^{n-1}$ the Hopf projection. 
Given an open cover by balls $\cu:=\{U_\alpha\}$ of $\pi(K\priv \cc)$,
choose an open set $U$ which completes $\cu$ to a cover of
$\P^{n-1}$. Consider a smooth partition of unity $\{\Phi_\alpha,\Phi\}$
associated to $(U_\alpha,U)$. Also choose a smooth non-negative function
$\theta$ on $\P^{n-1}$ which vanishes on $\pi(\cc)$ at order two and is
positive out of $\pi(\cc)$. We claim that if we can find functions
$h_\alpha:\pi^{-1}(U_\alpha)\lra \R$ such that $i\vec N\cdot h_\alpha<0$ on
$K\priv \cc\cap \pi^{-1}(U_\alpha)$
then the function $H=\theta \circ \pi \sum \Phi_\alpha\circ \pi\cdot
h_\alpha$ is of the required type. Actually, 
\begin{eqnarray}
i\vec N\cdot H &=& \theta \circ \pi\sum \Phi_\alpha\circ \pi\, i\vec
N\cdot h_\alpha <0 \text{ on } K\priv \cc\\
\text{and } dH &= & \left(\sum \Phi_\alpha\circ \pi. h_\alpha\right)
d(\theta \circ \pi)+\theta \circ \pi d\left(\sum \Phi_\alpha\circ
\pi. h_\alpha\right)=0 \text{ on } \cc.\label{eq12} 
\end{eqnarray}
The first equality holds because $\theta\circ \pi$ and
$\Phi_\alpha\circ \pi$ are constant along $i\vec N$. The second
because both $\theta \circ \pi$ and $d\theta \circ \pi$ vanish on $\cc$. 

Thus it only remains to build the cover of $\pi(K\priv \cc)\subset
\P^{n-1}$ by balls and the suitable functions on the solid tori which
projects to these balls by $\pi$. In the neighbourhood of $x\in
\pi(K\priv \cc)$, the compacity of $K$ allows to find a local section
$s:B_{\delta(x)}\lra S^{2n-1}$ of $\pi$ and constants $\eps(x)>0$ such 
that the intervals 
$$
I(\wdt x):=\left\{ \Phi_{i\vec N}^t(s(\wdt x))\; , \; t\in
]-\eps(x),\eps(x)[\right\},\hspace{.5cm} \wdt x\in B_{\delta(x)} 
$$ 
contain no points of $K$. These balls $B_{\delta(x)}$ provide a cover
of $\pi(K\priv \cc)$, for which we construct now the functions
$h_x$. For this, first identify $\pi^{-1}(B_{\delta(x)})$ with
$B_{\delta(x)}\times S^1$ by the map $\phi:B_{\delta(x)}\times
S^1\lra\pi^{-1}(B_{\delta(x)})$ defined by $\phi(\wdt x,t):=\Phi_{i\vec
  N}^{t+\eps(x)}(s(\wdt x))$. By $\phi$, $I_{\wdt x}$ is taken to
$]-2\eps(x),0[$ and $i\vec N$ to $\partial /\partial t$. In these
    coordinates, we can define the function $h_x$ on $B_{\delta(x)}\times
    S^1$ by :
$$
\left\{
\begin{array}{l}
h_x(\wdt x,t)=-t \hspace{.5cm} \text{for } t\in [0,1-2\eps(x)],\\
h_x \text{ is smooth on } B_{\delta(x)}\times S^1.
\end{array}
\right.
$$
Such a  function can obviously be defined, and it fits with the
requirement $i\vec N\cdot h_\alpha<0$ on $K\priv \cc\cap
\pi^{-1}(U_\alpha)$.\hfill $\square$
\begin{cor}[Perturbation procedure for smooth balls]
\label{pertsmooth} Let $M$ be a symplectic manifold, $B\subset M$ a
smooth symplectic closed ball and $U$ an open set of $M$. Assume that
$U\cap B=\emptyset$, denote $K:=\partial U\cap B$ and $\cc$ the Hopf
circles of $\partial B$ contained in $K$. Then there is a Hamiltonian
function $H$ on $M$ such that $\Phi_{X_H}^\eps(B)\cap
\overline{U}=\Phi_{X_H}^\eps(B)\cap \partial U=\cc$ for any small positive
real number $\eps$.    
\end{cor}
\noindent{\it Proof :} It is well-known that there exists a
symplectomorphism $\Phi$ between a neighbourhood $N(\partial B)$ of
$\partial B$ in $M$ and  a neighbourhood $N(S)$ of $S(r)$ in $\C^n$. This
map sends $K$ to $K'=\Phi(K)$, $\cc$ to $\cc'=\Phi(\cc)$ where $\cc'$
is exactly the set of Hopf circles of $S(r)$ contained in $K'$. Consider
a function $h$ on $S(r)$ associated to $(K'$, $\cc')$ as in lemma
\ref{sul}. Extend it to a smooth function of $\C^n$ with compact
support in $N(S)$. The function $H:=h\circ \Phi$, {\it a priori} defined
on $N(\partial B)$, can be extended to $M$ by setting $H=0$ outside
$N(\partial B)$. As explained above, the vector field $X_H$ points
inside $B$ on $K\priv \cc$ and vanishes on $\cc$. \cqfd
\indent This corollary is unfortunately only useful when dealing with smooth
balls. In our context, $B$ will 
usually be one of the (only regular) balls $B_i$ and $U$ the union of
the other balls $\cup_{j\neq i}B_j$. We will prove the counterpart of
this perturbation procedure for regular symplectic closed balls in
next section  (see proposition \ref{pertreg}).  \vspace{.3cm}\\
\textbf{Remark : a $\cc^0$ definition of the closed characteristics.} 
Sullivan's lemma suggests an alternative definition of a closed 
characteristic which makes no mention on the charcteristic distribution. 
Given an open set $U$ with smooth boundary and a point $x\in \partial U$, 
define 
$$
\begin{array}{l}
\ds C_x:=\overline{\cap_{\eps>0} C_{x,\eps}}\\
\ds C_{x,\eps}:=\left\{y\in\partial U\, | \, \exists H\in \cc^\infty(M,\R),\;
\begin{array}{l}
\Phi_{X_H}^t(x)\in U,\\
\text{and } \Phi_{X_H}^t(\partial U\priv B_\eps(y))\subset \overline{U} 
\end{array}
\hspace{0.5cm}\forall t>0\right\}.
\end{array}
$$
The set $C_x$ may be $\partial U$ itself for some points. However, Sullivan's 
lemma ensures that it is the characteristic leaf through $x$ when it is closed.
It would be interesting to know wether this definition gives a symplectic 
invariant when $\partial U$ is not smooth any more. For instance, do the $C_x$ 
and the Hopf circles coincide for a sufficiently regular but non-smooth symplectic
ball? Equivalently, one can ask if this definition gives non-trivial subsets of 
$\partial U$ when continuous Hamiltonians rather than smooth ones are considered
\cite{oh1}.
\subsection{Symplectic spheres in maximal packings.}\label{partsphere}
We show here that  any Hopf circle in the
intersection of two symplectic balls along their boundaries gives rise
to its own supporting sphere.
\begin{lemma}\label{sympsphere}
Let $B_1$, $B_2$ be regular symplectic closed balls of a symplectic
manifold  $(M,\om)$ with disjoint interior, of radii $r_1$, $r_2$ and
centers $O_1$, $O_2$. Through any common Hopf circle $C$ of $\partial B_1$
and $\partial B_2$ there is a topological 2-sphere $S_C$ of $M$ passing through $O_1$,
$O_2$, smooth except along $C$ and with symplectic area
$\pi(r_1^2+r_2^2)$. In dimension four, any two such spheres intersect precisely
in $O_1$ and $O_2$ with intersection number $1$ at each point. 
\end{lemma}
\noindent\underline{Proof :} Let $\phi_i:\overline{B(r_i)}\lra B_i$
the corresponding symplectic embeddings and
$C_i:=\phi_i^{-1}(C)$. Then $C_i$ is a Hopf circle of $\partial B(r_i)$
and bounds a holomorphic disc $\D_{C_i}\subset B(r_i)$. Its image by
$\phi_i$ is a symplectic disc $D_i$ in $B_i$ of area $\pi
r_i^2$ passing through $O_i$ and bounded by $C$ (a Hopf disc). 
When gluing $D_1$ with
$D_2$ along their common boundary, we get a topological
$2$-sphere $S_C$ which is smooth except on $C$ and has area
$\ca(S_C)=\pi r_1^2+\pi r_2^2$. By construction, two such spheres
$S_C$, $S_{C'}$ intersect only at $O_1$, $O_2$ and in dimension $4$,
their intersection 
numbers at $O_i$ is $1$ because they are the same as the intersection
numbers of $\D_{C_i}$, $\D_{C_i'}$ in $B^4(r_i)$. \hfill $\square$
\subsection{Smooth maximal packings : proof of theorem
  \ref{smooth}.}\label{partsmooth} 
\paragraph{Proof of theorem \ref{smooth}.a.}
Fix a smooth maximal symplectic packing of $\P^2$
 by two balls $\{B_1,B_2\}$ of radii $r_1,r_2$ checking the maximality 
condition $r_1^2+r^2_2=1$. Denote by $K:=\partial 
B_1 \cap \partial B_2$ and $\cc$ the set of Hopf circles of $\partial
B_1$ contained in $K$. Note that such a circle is also a Hopf circle
of $\partial B_2$ because $\partial B_1$ and $\partial B_2$ are
tangent along $K$. The perturbation procedure \ref{pertsmooth} shows
that after a possible Hamiltonian perturbation of the packing, we can
assume that $K=\cc$. 

Gromov's work shows that the balls of a smooth maximal packing of 
$\P^2$ by two balls cannot be disjoint (we also prove it for proposition 
\ref{extension} in the more general setting of regular mapping). 
It implies of course that no Hamiltonian perturbation of our packing can 
lead to disjoint balls. In view of the perturbation procedure \ref{pertsmooth},
$K$ must contain at least one Hopf circle of $\partial B_1$,
so $\cc$ is not empty. To prove that $\cc$ is exactly one Hopf circle,
we argue by contradiction and assume that $\cc$ contains two circles $C$ and $C'$. The
topological 2-spheres $S_C$ and $S_{C'}$ given by Lemma
\ref{sympsphere} intersect precisely at $O_1$, $O_2$ with
intersection numbers one at each point. On the other hand,
they both have symplectic area $\pi(r_1^2+r_2^2)=\pi$ so they are in
the homology class of a line $L$ in $\P^2$. We thus have
$2=S_{C_1}\cdot S_{C_2}=L\cdot L=1$ which is a contradiction.\hfill $\square$
\paragraph{Proof of theorem \ref{smooth}.b,c.} 
Let $B_1,B_2,B_3$ be the closed balls of our smooth maximal
packing. We can assume that their radii check 
$r_1\geq r_2\geq r_3$. Then the maximality condition is :
$$
r_1^2+r_2^2=r_1^2+r_3^2=1. 
$$ 
Define as above $K_{ij}:=B_i\cap B_j=\partial B_i\cap \partial B_j$
and the sets $\cc_{ij}$ consisting 
of the Hopf circles of $K_{ij}$.
The first point is to prove that we can get rid of $K_{ij}\priv
\cc_{ij}$ by a Hamiltonian perturbation. Unlike the two-balls
situation, it is not obvious 
at first glance because a Hopf circle of $\partial B_1$ could {\it a
  priori} be covered by $\partial B_2\cup \partial B_3$ without being
in any $K_{ij}$. But notice that the $K_{ij}$ are pairwise
disjoint compact sets (and remain so after perturbation) because
the intersection of two of the $K_{ij}$ is precisely $B_1\cap B_2\cap
B_3$. This intersection must be empty because any point of it would be
a singular point of at least one of the balls.
By connexity, we deduce that $K_{12}\cup K_{13} $ contains no Hopf
circles away from those in $\cc_{12}$ and $\cc_{13}$. After a 
Hamiltonian perturbation of $B_1$ according to corollary
\ref{pertsmooth}, we can assume that $K_{i1}=\cc_{i1}$. After a second
Hamiltonian perturbation, this time of $B_2$, we can also assume that
$K_{23}=\cc_{23}$.\\
\textit{Proof of theorem \ref{smooth}.b ($r_1=r_2=r_3=1/\sqrt 2$).}
In this case, any two balls of $B_1$, $B_2$, $B_3$ form a maximal
symplectic packing of $\P^2$ by two balls. As such, and in view of
theorem \ref{smooth}.a, $K_{ij}=\cc_{ij}$ contains exactly one 
circle.
\vspace{.2cm}\\
\textit{Proof of theorem \ref{smooth}.c ($r_1>r_2=r_3$).} In this case,
$(B_1,B_2)$ and $(B_1,B_3)$ form maximal symplectic packings of $\P^2$
by two balls. As behind we conclude that $\cc_{12}$ and $\cc_{13}$
contain exactly one circle. Finally, if $\cc_{23}$ were to contain a
circle, Lemma \ref{sympsphere} would associate to it a toplogical
$2$-sphere of symplectic area $\pi(r_2^2+r_3^2)\in ]0,\pi[$. It is
clearly impossible, so $\cc_{23}$ is empty.\hfill $\square$
\begin{rk}\label{pratique}
This last argument shows that any two closed symplectic balls of radii
$r_1,r_2$ in $\P^2$ with disjoint interiors can intersect along a full Hopf
circle of the boundary of one of them (henceforth of both) only if
$r_1^2+r_2^2\in \N$. 
\end{rk}
\paragraph{Proof of theorem \ref{smooth}.d}
 Suppose by contradiction that $n$ closed balls $B_1,\dots,
B_n$ of the same radius $r$ constitute a
maximal packing of $\P^2$ ($n\geq 4$).

We know from \cite{mcpo} and \cite{biran2}
that for $n=4$ and $n\geq 9$, such a packing should fill the
space. Then any point $p$ of the boundary (in $\partial B_1$) of
$\partial B_1\cap \partial B_i$ belongs to a third ball $B_j$. So our
packing by $B_1,\cdots,B_n$ is not smooth as already noticed in last
paragraph. Moreover, the radius of the balls is less than $\sqrt{2/5}<\sqrt{1/2}$ 
for $n=5,6,7,8$. As before, we are in a non-removable intersection situation : 
the ball $B_1$ cannot be disjointed from the union of the other ones $\cup_{j\neq
  1}B_j$ (see proposition \ref{extension}). In view of corollary
\ref{pertsmooth}, it means that there is a Hopf 
circle $C$  of $\partial B_1$ which is also covered by $\cup_{j\neq 1}
\partial B_j$. Assume to fix the  notations that $C\cap \partial
B_2\neq \emptyset$. Since $r^2<1/2$, remark \ref{pratique} shows that
$C$ cannot be entirely contained in $\partial 
B_2$. So there is a point $p$ in the boundary of $\partial B_2$ in
$C$. Such a point also belongs to another ball, so it is an
intersection point of $\partial B_1,\partial B_2$ and $\partial B_i$
for $i\neq 1,2$. Our packing cannot be smooth.\hfill $\square$                      
\section{Regularizations and extensions of regular embeddings.}\label{partregularization}
This purely technical section is aimed at extending the perturbation
procedure described in corollary \ref{pertsmooth} to non-smooth
balls. We first explain how to extend slightly a
symplectic open ball in an open manifold $M$. It should be
noticed that the following proposition applies to regular embeddings
of a ball.
\begin{prop}\label{extension}
Let $\phi:B(r)\hra M$ be a symplectic embedding of an euclidean
ball. Assume $\phi$ extends smoothly to an embedding of 
$\overline{B(r)}\priv K$ into $M$ where $K$ is a closed set of
$\partial B(r)$. If $K$ contains no Hopf circle of $\partial B(r)$
then for any open neighbourhood $V$ of $\overline{\im \phi}$, there is
a symplectic embedding $\wdt \phi_\eps:B(r+\eps)\hra V$ for $\eps$
small enough.  

If $K=\emptyset$ then $\wdt \phi_\eps$ can be choosen to be an
extension of $\phi$ {\it i.e.} $\wdt \phi_{\eps\,|B(r)}=\phi$.  
\end{prop}
The non-squeezing theorem shows that this proposition is sharp in the
sense that $K$ actually has to be assumed not to contain any Hopf
circle. \\  
\noindent {\it Proof :} When the singular locus is empty
($K=\emptyset$), this proposition is very classical. We show however
its brief proof so that it serves as a basis for the non-smooth
case. In this situation, the image $\phi(\partial 
  B(r))$ is a smooth hypersurface of $M$. Denoting by $S:=\partial B(r)$,
the regular neighbourhood theorem gives a symplectomorphism
$\psi:N_\eps(S)\lra N(\phi(S))$ from an $\eps$-neighbourhood of the 
standard sphere of radius $r$ in $\C^n$ to a neighbourhood of $\phi(S)$
 in $M$ (see \cite{weinstein} or \cite{mcsa}, p.101). 
Morevover, this map can be chosen in such a way that
$\Psi_{|S}=\phi_{|S}$ and $\psi'(x)\vec N(x)=\phi'(x)\vec N(x)$ for
any point $x\in S$ ($\vec N(x)$ denotes again $x/\Vert x\Vert$). The
map $\phi_\eps:B(r+\eps)\lra M$ defined by 
$$
\phi_\eps(x):=\left\{
\begin{array}{l}
\phi(x) \text{ for } x\in B\\
\Psi(x) \text{ for } x\in N_\eps(S)\priv B  
\end{array}\right.
$$  
is $\cc^1$-smooth, and hence a symplectomorphism from the ball
of radius $r+\eps$ inside $M$.   

For singular hypersurfaces, the standard neighbourhood theorem is not 
  valid anymore, hence the above proof has no straightforward
  generalization. However it can
  be easily adapted by using a regularization trick explained in
  Lemma \ref{regularization} below. It states that under the extension
  condition on $\phi$ of proposition \ref{extension}, there exists a
  ``regularization'' $\wdt \phi: \overline{B(r)}\hra M$ which is a
  symplectic embedding of the {\it closed} ball inside the prescribed
  neighbourhood $V$ of $\phi(\overline{B(r)}\priv K)$. Now the
  extension $\wdt \phi_\eps$ of $\wdt \phi$ constructed in the case
  $K=\emptyset$ gives the desired map.  \cqfd
The result we have used in this proof is a particular
  case of the following lemma for $\cc=\emptyset$. This extended
  version will nevertheless be usefull in order to generalize
 corollary \ref{pertsmooth}.
\begin{lemma}\label{regularization}
Let $\phi:\overline{B(r)}\priv K\hra M$ be a symplectic embedding of a
closed ball minus a singular compact set $K\subset \partial B(r)$
inside an open manifold $M$. Denote by $\cc$ the set of Hopf
circles of $\partial B(r)$ contained in $K$. There exists a
symplectic embedding $\wdt \phi$ of  $\overline{B(r)}\priv \cc$ inside
any prescribed neighbourhood $V$ of $\im
\phi=\phi(\overline{B(r)}\priv K)$.  
\end{lemma}
\noindent {\it Proof :} Consider an open set $U\Subset \partial
B(r)\priv K$ which contains at least one point of each Hopf circle of
$\partial B(r)\priv \cc$. As in the previous proof, $\phi$ can be
extended as a symplectic embedding to a shell
$U_\eps:=\{[x,(1+\eps(x))x[,\; x\in U\}$, where $\eps$ is a small
positive function on $U$. Provided $\eps$ is sufficiently small, $\phi$
actually sends $B_\eps':=B(r)\priv K\cup U_\eps$
into $V$. Lemma \ref{sul} shows that there exist a smooth function $h:\partial
B(r)\lra \R$ such that $i\vec N\cdot h<0$ on $\partial B(r)\priv
(U\cup \cc)$ with $dh=0$ on $\cc$. When extending this function to
$\C^n$, we get  a Hamiltonian function whose flow has the property
that for $t$ small enough, the set $\Phi_{\vec X_h}^t(\overline{B(r)}\priv
\cc)$ is contained inside $B_\eps'$. The map  
$$
\wdt \phi:=\phi\circ \Phi_{\vec X_h}^t:\overline{B(r)}\priv \cc\hra V 
$$    
gives the desired symplectic embedding of $\overline{B(r)}\priv \cc$
inside $V$.\cqfd
Before stating the anounced generalization of the perturbation procedure
\ref{pertsmooth} to regular balls, we need to broaden slightly the
notion of Hamiltonian perturbation. We say that a path $\{\phi_t,B_t\}$
of regular symplectic balls in $M$ is a {\it regular Hamiltonian
  deformation} if there exist smooth functions
$H_t$ defined in the interior of $B_t$ such that
$\phi_t=\Phi_{X_{H_t}}^t\circ \phi$. The point is that we do not
impose to these functions to be defined in all of $M$, so they may
have singularities on the boundary. We however demand the whole path
to be made of regular symplectic balls. Now, a regular Hamiltonian
perturbation of a packing is a regular Hamiltonian perturbation of
each of its balls such that the configuration of the balls at any time
$t$ remains a packing. To illustrate this definition, consider the
special case of lemma \ref{regularization} when $M$ is an open
manifold with  boundary $\partial M$ and $\phi:(B(r),K)\hra
(M,\partial M)$ is a regular symplectic ball with a part $K$ of its 
boundary sent in $\partial M$. Then the map $\wdt \phi$ is in fact a regular
Hamiltonian deformation of $\phi$ associated to the function $h\circ
\phi^{-1}$. Moreover, $\wdt \phi_{|\cc}=\phi_{|\cc}$. This obvious
remark is exactly the content of next proposition.  
\begin{prop}[Perturbation procedure for regular balls]\label{pertreg}
Let $U$ be an open set of a symplectic manifold $M$ and $B$ a regular
symplectic ball in $M\priv U$. Denote $K:=\partial B\priv
\partial U$ and $\cc$ the set of Hopf circles of $\partial B$
contained in $K$. Then there is a regular Hamiltonian perturbation
$(B_t)_{t<\eps}$ of $B$ with $B_t\cap \partial U=\cc$ for any positive
$t$.
\end{prop}
As a corollary of proposition \ref{extension}, we get a non-removable
intersection property for regular maximal symplectic packings. 
No ball of such a packing is completely disjoint from the 
other balls. 
\begin{cor}[Non-removable intersection]\label{nri} Let 
$B_i:=\phi_i(\overline{B(r_i)})$ be a regular maximal symplectic packing of
$M$. Then
$$
\forall i\in [1,k], \hspace{.5cm} B_i\cap(\cup_{j\neq i}B_j)\neq \emptyset.
$$ 
This intersection even contains at least one Hopf circle of $\partial
B_i$.
\end{cor} 
\noindent {\it Proof :} Call $\phi_i$ the symplectic embeddings of $B(r_i)$
in $M$ corresponding to $B_i$. Suppose by contradiction that 
$K:=B_i\cap(\cup_{j\neq i}B_j)$ contains no Hopf circle of $\partial
B_i$. The map $\phi_i$ is a regular symplectic embedding of
$\overline{B(r_i)}\priv K$ inside the open symplectic manifold
$M'=M\priv \cup_{j\neq i}B_j$. By proposition \ref{extension},
$\phi_i$ can be ``extended'' to a symplectic 
embedding  $\wdt \phi_i:B(r_i+\eps)\hra M'$. The maps 
$\phi_1,\dots,\wdt \phi_i,\dots, \phi_k$ thus provide a symplectic
packing of $M$ by $k$ balls of radii $(r_1,\dots,
r_i+\eps,\dots,r_k)$. It is clearly in contradiction with the
definition of a maximal symplectic packing.  \hfill  $\square$ 
\section{Non-smooth symplectic packings.}\label{partreg}
We turn to the problem of identifying non-smooth (but regular) maximal 
symplectic packings. We first prove the existence of supporting surfaces of such
packings  (theorem \ref{regular}). We then remark that passing from theorem 
\ref{regular} to the precise statement of theorem \ref{theend} is a matter
of being able to perturb the topological supporting surface to a smooth symplectic 
surface (lemma \ref{dream}). We finally show that the smoothening is possible
under the hypothesis of theorem \ref{theend}, thus proving it. 
\subsection{Existence of a supporting surface : proof of theorem \ref{regular}.}
The idea of the proof is very simple. The non-removable intersection
property of the packing implies the existence of a Hopf circle $C_1$
of $\partial B_1$ contained in the union of the other balls. Let us
denote $\cc_{1i}$ the set of circles of $B_i$ intersecting $C_1$ along an open
set ($\cc_{1i}$ may be empty or contain several circles). Denote  $S_0$
the Hopf disc in $B_1$ bounded by $C_1$ and $S_1$ the surface obtained
by gluing the Hopf discs corresponding to  $\cc_{1i}$  to
$S_0$. Obviously no point of $C_1$ is a boundary point of $S_1$. Now
there are two possibilities. Either the $\cc_{1i}$ are covered by the
balls $B_j$ for $j\neq i$ or not. In the first case, we can glue to
$S_1$ the  Hopf discs corresponding to the circles of $\cc_{1ij}$
($\cc_{1ij}$ is the set of circles of $ \partial B_j$ intersecting
a circle of $\cc_{1i}$ in an open set). We obtain a surface $S_2$ with boundary
points neither in $C_1$ nor in $C_{1i}$, which we can use in order to
iterate the construction. In the latter, we can get rid of one of the
circles of $C_{1i}$ from the intersections between the balls by a small
Hamiltonian perturbation. The effect of this transformation is that
$C_1$ is not covered by the other balls any more. It actually means
that $C_1$ was not relevant for our purpose, but there clearly exists
another circle of $\partial B_1$ to which we can apply the previous
procedure. The reason for which this iteration process stops and
produces a closed surface is the finiteness condition on the singularities. 
In order to make rigorous and clear the preceding iterative process,
we encode the situation into a graph. 

Let $B_1,\dots,B_k$ be the balls of a regular maximal symplectic packing 
of $M$. We already know how to associate a supporting topological
sphere to any common Hopf circle of two balls of the packing (see
lemma \ref{sympsphere}). We thus consider in the following only those
Hopf circles not concerned by this basic construction. Define 
$$
\begin{array}{l}
\sc_i:=\{x\in \partial B_i\; |\; \exists j\, ,\; C_x\subset \partial B_j\},\\
\cc_i:=\{x\in \partial B_i\; |\; C_x\subset \bigcup_{j\neq i} \partial
B_j\;,\; C_x\not\subset \partial B_j\;\; \forall j\neq i\},
\end{array}
$$ 
where $C_x$ denotes the possibly singular Hopf circle of $\partial B_i$ 
passing through $x$. Clearly, if one of the $\cc_i$ is empty, then
corollaries \ref{nri} and \ref{sympsphere} show that there must be a
supporting $2$-sphere of the packing passing through $B_i$. As before,
our aim is to explain to what $\cc_i$ can be reduced after Hamiltonian
perturbation of the packing. Notice that by definition, each circle  
of $\cc_i$ contains at least one triple intersection point between the
balls of the packing, and there are only finitely many such points
because of the regularity condition on the packing. Since any two
Hopf circles of a given ball are disjoint, each $\cc_i$ thus  
contains only a finite number of Hopf circle of $\partial B_i$.  
Consider henceforth the finite graph $G$ whose vertices are the Hopf
circles contained in one of  the $\cc_i$,  and the edges are the pairs
of such circles which share an open 
arc. Also colour the vertices black when they represent a circle $C\in
\cc_i$ which is also contained in $\cup_{j\neq i} \cc_j$ and red otherwise.

A red vertex is a Hopf circle $C$ of $\partial B_i$ which is covered by 
the other balls but not by the union of the $\cc_j$ for $j\neq
i$. Since each $\cc_j$ is compact, there must be an open arc $I\subset
C$ which is a piece of a Hopf circle of $\partial B_j$, not covered
itself by the other balls of the packing. The perturbation procedure
then allows to produce a new packing $(\wdt B_1,\dots,\wdt B_k)$ very
close to the original one, with intersection between the
balls unchanged except that $\wdt B_i\cap \wdt B_j=B_i\cap B_j\priv I$.
The graph $\wdt G$ associated to the perturbed maximal packing is thus
 a subgraph of $G$ obtained by erasing the vertex $C$ together with all 
its adjacent edges, and turning all its neighbours in $G$ to red. In 
particular, $\wdt G$ has one vertex less than $G$. This process can be 
iterated as long as there is a red vertex in the graph. Because the initial 
graph is finite, there must be a stabilization after a finite number of 
steps. We conclude that some Hamiltonian perturbation of the packing leads
to a graph which is only black. We will now suppose that $G$ itself
has only black vertices. Applying the perturbation procedure once
more, we can arrange so that 
$$
B_i\cap \cup_{j\neq i} B_j=\cc_i\cup\sc_i.
$$
 
We claim that each connected component $\hat G$ of $G$ corresponds
to a supporting surface. Actually, let $S$ be the union of the Hopf
discs corresponding to the Hopf circles of $\hat G$. It is obviously
a connected space, covered by the closed balls of the packing. We need
however a brief discussion of the regularity of this space before we
can assure it is actually a topological surface. 
Inside the balls first, $S$ is an immersed symplectic surface, whose
only self-intersection points are positive and located at the origin
of the balls. Consider then a point $x$ of $S\cap \partial B_i$ which is
not a singular point of any ball. In particular it is not a triple
intersection point of the packing, so $S$ cannot be made of more than
two discs in a neighbourhood. Moreover $x$ belongs to a Hopf circle
of $\cc_i$. The fact that this Hopf circle is black coloured implies
precisely that $x$ is in the boundary of at least two Hopf discs.
It follows that $S$ is  locally made of exactly two smooth Hopf
discs glued along a interval around $x$. It is easily checked that $S$
is even locally diffeomorphic to a cylinder $\{y=|x|\}\times \R$ at these
points. We thus conclude that $S$ is as announced a closed topological
surface with finitely many possible singular points located at the
singular points of the packing.\hfill $\square$
\subsection{Refinement of theorem \ref{regular} when the supporting surfaces are 
smooth.}
In the concrete case of five balls in $\P^2$, we explain now how to sharpen
theorem \ref{regular} and get theorem \ref{theend}.   
We consider for the remaining of this paragraph a regular maximal
symplectic packing of $\P^2$ by five equal balls $\phi_1,\dots,\phi_5$ (or
$B_1,\dots,B_5$) - of radii $\sqrt{2/5}$ - and a supporting surface
$S$ of the packing. Then 
$S$ is made of several Hopf discs, each of which is of area $2/5\pi$ because
its boundary is a Hopf circle of the boundary of a ball.  Since the
Fubini-Study form is integral and $S$ is a closed surface, the number
of discs is a multiple $5k$ of five, the symplectic area of $S$ is
$2k\pi$, and its homology class is $2k[L]$ ($L$ is a line in $\P^2$).  
\begin{lemma}\label{dream}
Assume there is a smooth symplectic immersion $\cc^0$-close to $S$, whose
only self-intersections are at the origins of the balls. 
Then $S$ is of area $2\pi$, made of one Hopf disc in each ball of the packing, 
and it is the unique supporting surface of the packing. 
\end{lemma} 
\noindent{\it Proof :} Denote by $k_i$ the number of Hopf discs of
$S\cap B_i$, so that $k_1+k_2+\dots+k_5=5k$. The self-intersection of
$S$ in a  neighbourhood of the origin $O_i$ of $B_i$ is then given by
the formula :
$$
\delta_i=\frac{k_i(k_i-1)}{2}.
$$
The assumption on the smoothening of $S$ means that one can find a 
symplectically immersed surface $\wdt S$ homologous to $S$ (hence in the 
homology class $2k[L]$) with positive self-intersections. They are located 
in small neighbourhoods of the $O_i$ and are the same as those of $S\cap B_i$.
The total self-intersection number of $\wdt S$ is thus
$$
\delta:=\sum_{i=1}^5 \delta_i=\sum_{i=1}^5 \frac{k_i(k_i-1)}{2}.
$$
Taking into account that $\sum k_i=5k$, we easily get that 
\begin{equation}\label{five1}
\delta\geq \frac{5k(k-1)}{2}.
\end{equation}
The positivity of the self-intersections together with the fact that 
 $\wdt S$ is symplectic imply
  that $\wdt S$ is actually a $J$-holomorphic curve for a good almost
  complex structure on $\P^2$ see \cite{aula}. It must therefore 
verify the adjunction inequality, which gives in our present situation :
\begin{equation}\label{five2}
\delta\leq \frac{(2k-1)(2k-2)}{2}=(2k-1)(k-1).
\end{equation}
It follows from (\ref{five1}) and (\ref{five2}) that $k=1$, and
$\delta=0$. We thus conclude that $S$ was made of five Hopf discs, one
in each ball. Finally, we argue by contradiction to prove that $S$ is
the only supporting surface of the packing. Assume that there is a
supporting surface $S'$ distinct of $S$, made of $k_i'$ Hopf discs in
each balls, and  of total symplectic area $2k'\pi$ (so that $\sum
k_i=5k'$). Since the only intersection points between $S$ and $S'$ are
the center of the balls, we get :
$$
4k'=2L\cdot 2k'L=S\cdot S'=\sum_{i=1}^5 S\cdot_{O_i}S'=\sum_{i=1}^5
k_i\cdot 1=5k'.
$$  
This is the desired contradiction.\cqfd 
Note that the previous computation can be made for seven ({\it resp.} eight) balls. 
When they are symplectically smoothable in the previous sense, the 
supporting surfaces are at most seven ({\it resp.} eight), all of area $3\pi$ 
({\it resp.} $6\pi$). Each one intersects six of the balls 
through one Hopf disc and the last one through two Hopf dics ({\it resp.} intersects 
seven of the balls through two Hopf discs and the last one through three). 
\subsection{Smoothening of the supporting surfaces for packings of simple type.}
In this paragraph, we show that the smoothening required by lemma \ref{dream} 
can be achieved when precise conditions on the singularities of the packings are 
given. Recall 
that the singularities of $S$ are the union of a finite set of singularities 
of the packing and segments of the characteristic foliations joining precisely
these points. The first step is to deal with the generic singular points of $S$ .
\begin{lemma}\label{remseg}
Let $S$ be a symplectic surface singular along a segment $\Gamma$, a neighbourhood of 
which is diffeomorphic to 
\begin{equation}\label{nicecyl}
\{y=\alpha(z)|x|\}\subset \R^2(x,y)\times ]0,1[(z),
\end{equation} 
where $\alpha(z)$ is a continuous function which vanishes at $0$ and $1$. 
Then $S$ can be smoothened to a symplectic surface by a $\cc^0$-perturbation.
\end{lemma}
\noindent{\it Proof :} Consider such a segment of singularities $\Gamma$ of $S$. 
Using Moser's argument, a neighbourhood of $\Gamma$ can be presented symplectically 
as the cylinder $V_\eps:=\{|x_1|<1+\eps\;, \; |y_1|<\eps\; ,\; |z_2|<\eps\}\subset 
\C^2(z_1=x_1+iy_1,z_2)$ in such a way that $\Gamma$ corresponds to $[-1,1]\times 
\{0\}$ and $S$ corresponds to a union of two symplectic surfaces $S_1$, $S_2$ with 
common boundary $\Gamma$. Since $S_1$ and $S_2$ are symplectic, the identification 
can be done in order to achieve also 
$T_qS_1=\span_\R(\partial/\partial x_1,\partial/\partial y_1)$ and 
$T_qS_2=\span_\R(\partial/\partial x_1,-\partial/\partial y_1-u(x_1))$ 
($q\in \Gamma$) where $u(x_1)$ is a $\cc^1$-smooth vector field along $\Gamma$ with 
values in $\span_\C(\partial /\partial z_2)$ which vanishes for $x_1\notin ]-1,1[$. 
Now if the neighbourhood of $\Gamma$ is small enough, $S_1$ can be 
straightened  to the strip $A=\{|x_1|<1+\eps\, ,\, 0\leq y_1<\eps\, ,\, z_2=0\}$ 
by a map $h$ 
which is $\cc^1$-close to the identity (and even tangent to the identity along $\Gamma$). 
This map may distort $\om$ but by no more than an $\eps$-factor. We produce our 
symplectic smoothening of $S$ by cutting $S_1$ and replacing it by a ``very symplectic''
surface $\Sigma\subset V_\eps$ which interpolates smoothly between 
$\big(\Gamma,\Span_\R(\frac{\partial}{\partial x_1},
\frac{\partial}{\partial y_1}+u(x_1))\big)$ and $A$.

To this purpose, we consider a smooth ``profile'' of maps $\phi_{\vec v}:[0,\eps[\lra \C$
parameterized by vectors $\vec v$ of $\C$, such that $\phi_{\vec v}(0)=0$, 
$\dot\phi_{\vec v}(0)=\vec v$, $\phi_{\vec v}\equiv 0$ on $[\eps/2,\eps[$ and 
$\phi_0\equiv 0$. Up to shrinking these maps, we can also arrange so that 
$\partial \phi_{\vec v}/\partial \vec v$ and 
$\dot\phi_{\vec v}$ are very small on big compact sets (for $\vec v$). Then the surface 
$$
\Sigma:=\{\big(x_1,y_1,\phi_{u(x_1)}(y_1)\big),\; |x_1|<1+\eps\; , \; 0\leq y_1<\eps,\}
$$ 
obviously interpolates smoothly between $(\Gamma,\partial/\partial y_1+u(x_1))$ and 
$A$. Observe moreover that the tangent vectors to $\Sigma$
$$
v_1:=\frac{\partial}{\partial x_1}+\frac{\partial \phi_{\vec v}}{\partial \vec v}\cdot\frac
{\partial u}{\partial x_1}\; \hspace{1cm} \text{and} \hspace{1cm} 
v_2:=\frac{\partial}{\partial y_1}+\frac{\partial \phi_{\vec v}}{\partial y_1}
$$
are small perturbations of $\partial/\partial x_1$ and $\partial/\partial y_1$ 
respectively provided that $\phi$ is sufficiently small in $\cc^1$-norm. 
The tangent planes to $\Sigma$ are thus far from being Lagrangian\cqfd   
Hence all the difficulty of our desingularization problem concentrates at the 
(unavoidable) intersection of $S$ with singular points of the packing. To understand 
the situation in the greatest generality, we need local models for the
singularities that may arise in regular packings. We do not pretend to find them 
in this paper. Instead, and rather as an illustration, we focus on the very special 
type of singularities which appear in the examples we constructed in section 
\ref{examples}.
\begin{definition}\label{verysimple}
We say that a boundary singularity $p=\phi(q)$ of a symplectic ball $\{\phi,B\}$ is 
simple if $\phi$ is continuously differentiable at $q$, with non-vanishing derivative
in any directions but the charachteristic one. 
A regular packing is said to be of simple type if all the singularities are simple,
and if there exist no intersection point between any four balls (we only allow 
triple intersection points).
\end{definition}
Before proving theorem \ref{theend}, let us discuss the meaning of this definition
in our context. Consider a regular symplectic ball $\{\phi,B\}$ of a symplectic manifold 
and a simple singularity $p$ of $B$. Denote by $C_p$ and $D_p$ the Hopf circle and disc 
passing through $p$. For $q=\phi(\eta)\in C_p$, note also $T^c(q):=\phi'(\eta)\cdot 
T_\eta^\C\partial B$ and $\pi(q):=T^c(q)^{\perp \om}$. Then $D_p$ is tangent to $\pi(q)$
along $C_p$ and from the definition of a simple singularity, both symplectic plane 
distributions $T^c(q)$ and $\pi(q)$ have a well-defined limit $T^c(p)$ and $\pi(p)$ 
when $q$ goes to $p$. In particular, $D_p$ has a tangent plane $\pi(p)$ at $p$. 
Although $T_qC_p$ may not have limit at $p$, $C_p$ is tangent to $\pi(p)$ at $p$. 
Observing that the tangent plane to $\partial B$ at $q\in C_p$ is 
$T^c(q)\oplus T_qC_p$ we see that $\partial B\cap U$ is $\cc^0$-close to the cylinder 
$\pi(C_p)\times T^c(p)$ (and $B\cap U$ is close to $D_p\times T^c(p)$).
This remark immediately yields the following lemma :
\begin{lemma}\label{nospiral}
Let $B$ be a symplectic ball with a simple singularity at $p\in \partial B$. With the 
notations above, there is a neigbourhood $U$ of $p$ such that the linear projection 
$\pi:\overline{D_p}\cap U\lra \pi_p$ along $T^c(p)$ is an injective map. 
\end{lemma}
 This lemma prevents $D_p$ from spiraling too much above its tangent plane, 
creating problematic cone singularities. The proof of proposition 
\ref{theend} is achieved in two steps. First, we smoothen the supporting surface 
at each (simple type) singularity of the packing. Being cautious enough in the first 
step allows to use lemma \ref{remseg} to get rid of the remaining closed arcs of 
singularities.\vspace{.2cm}\\
{\it Proof of proposition \ref{theend} :} Consider a maximal symplectic packing 
$\{B_1,\dots,B_k\}\subset \P^2$ of simple type, and one of its supporting surfaces 
$S$. Denote $\{p_1,\dots,p_n\}$ the singularities of the packing which belongs to 
$S$. Then the singularities of $S$ consist of the points $p_i$ themselves, together 
with open segments of the characteristic foliation linking the $p_i$. The 
differentiable model of the singularities along these segments is locally 
$\{y=|x|\}\times \R\subset \R^3$. As announced, we are going to smoothen $S$ around
the $p_i$ in such a way that the remaining segments of singularities have the 
global form (\ref{nicecyl}). Consider henceforth one of the singular 
point $p:=p_i$. Recall that $p$ can only be an intersection between two or three 
balls. 
Assume first that $p$ lies in the intersection between two balls only, say $B_1$ and 
$B_2$. Let us distinguish between two cases.

{\it Case a : $\pi_1(p)=\pi_2(p)$, $T^c_1(p)=T_2^c(p)$.} All indices refer to the 
balls of the packing (for instance $\pi_1(p)$ is the plane defined above for the 
ball $B_1$). Consider local symplectic coordinates taking $p$ to the 
origin in $\C^2$ and $\pi(p)$
to $\{z_2=0\}$. Inside a small bidisc $Q_\eps:=\{|z_1|<\eps\}\times \{|z_2|<\eps\}$,
the surface $S$ together with its tangent planes are very close to $\{z_2=0\}$. The 
projection $\pi:S\cap Q_\eps\lra \{z_2=0\}$ is therefore a covering map, which must
be injective from lemma \ref{nospiral}. The intersection of $S$
with $\partial Q_\eps$ is therefore the graph over $\partial \D_\eps$ of a piecewise 
$\cc^1$-smooth complex valued function with small $\cc^1$-norm. 
This function can obviously be extended to $\D_\eps$ in such a way
that its graph $\Sigma\subset Q_\eps$ is tangent to $S$ on $\partial Q_\eps$ 
(except at the singularities of $S\cap \partial Q_\eps$), that its singularities 
are located on segments inside $\{\eps/2<|z_1|<\eps\}$ and have the form of 
lemma \ref{remseg}. Moreover, this extension can be choosen $\cc^1$-small,
so that $\Sigma$ remains a symplectic surface. Cutting $S\cap Q_\eps$ and replacing 
it by $\Sigma$ thus gives the desired smoothening of $S$ at $p$.

\begin{figure}[h!]
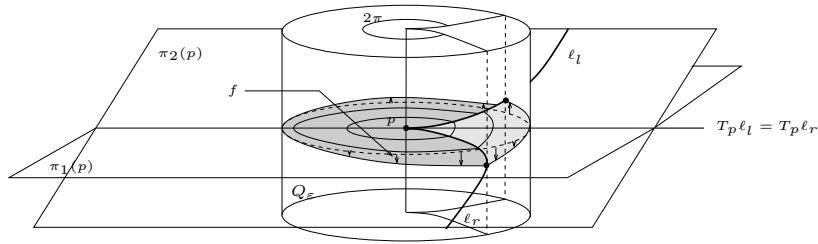

\begin{center} 
\input fun3.pstex_t
\caption{Singularity of type $\pi_1(p)=\pi_2(p)$ and $C_p$ is not smooth.}
\label{fun}
\end{center}
\end{figure}

{\it Case b : $\pi_1(p)\neq \pi_2(p)$.} Consider a small neighbourhood $U$ 
of $p$ such that the common Hopf circle $C_p:=C_{p1}=C_{p2}$ of $\partial B_1$, 
$\partial B_2$ passing through $p$ is the union of two smooth (open) arcs $\ell_l$
and $\ell_r$ meeting at $p$. Assume also that $\pi_1(q)\neq \pi_2(q)$ for all 
$q\in \ell_l\cup\ell_r$. Inside $U$ we have $T_qC_p=\pi_1(q)\cap \pi_2(q)$
so that $\ell_l$ and $\ell_r$ have common tangency $\pi_1(p)\cap \pi_2(p)$ 
at $p$. Since $D_i(p)$ has limit tangent plane $\pi_i(p)$ at $p$, 
if $\ell_l\cup \ell_r$ is $\cc^1$-smooth then the  model of the singularity 
of $S$ at $p$ is the same as at any generic point of $C_p$, and we can forget it.
Else $\ell_l\cup\ell_r$ is a $\cc^1$-cusp, meaning that its projection to 
$\pi_1(p)\cap \pi_2(p)$ is a half-line (see figure \ref{fun}). Then $C_p\times
 T_1^c(p)$ ($=C_p\times T^c_2(p)$) separates $U$ in two cylinders, a big one (with
aperture $2\pi$ at $p$) and a small one, in the whereabouts of which each ball is 
contained. To fix the ideas, suppose that $B_1$ is the ``big ball'' in $U$, and 
consider symplectic coordinates in $U$ such that $p=0$, $\pi_1(p)=\{z_2=x_2+iy_2=0\}$
and $\pi_1(p)\cap \pi_2(p)=\Span_\R(\partial /\partial y_1)$. Note that $\pi_2(p)$ is 
transverse to $\{z_1=0\}$ because it is symplectic and contains $\partial /\partial y_2$.
So the projection of $S$ on $\pi_1(p)$ along $T^c_1(p)$ inside a small bidisc $Q_\eps$ 
is a covering map, injective by lemma \ref{nospiral}.
 The intersection $\partial Q_\eps\cap S$ is henceforth the graph of 
a piecewise $\cc^1$-smooth complex valued function over $\partial \D_\eps$. It is 
easy to see that this map has bounded derivatives and small $\cc^0$-norm. The same 
procedure  as in case {\it a.} above produces the smoothening.  

\begin{figure}[h!]
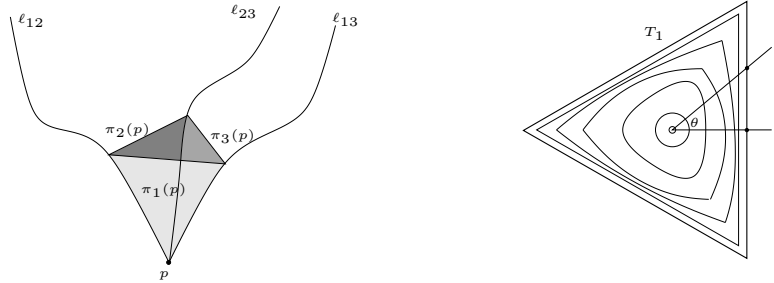

\begin{center} 
\input ouf.pstex_t
\caption{Singularity of $S$ and extrapolation between $T_1$ and $0$.}
\label{ouf}
\end{center}
\end{figure}

Assume now that $p$ lies in the boundary of three balls $B_1$, $B_2$, $B_3$. 
If looking inside
a sufficiently small neighbourhood $U$ of $p$, each intersection $B_i\cap B_j\cap U$ is
a smooth open arc $\ell_{ij}$ ending at $p$. Reasoning as in cases {\it a.} and {\it b.} 
it is easy to smoothen $S$ at $p$ when $\pi_1(p)=\pi_2(p)=\pi_3(p)$ or
$\pi_1(p)=\pi_2(p)\neq\pi_3(p)$. Thus it only remains to investigate the situation
of three different tangent planes for $D_i(p)$ at $p$ (see figure \ref{ouf}). 
The curves $\ell_{ij}$ 
have then well-defined tangencies at $p$ :
$$
T_p\ell_{ij}=\pi_i(p)\cap \pi_j(p).
$$  
Consider a parameterization of the $\ell_{ij}$ by smooth maps 
$\gamma_{ij}:[0,\eps[\lra U\subset \C^2$ with $\gamma_{ij}(0)=p$.
Notice that $\om(\dot\gamma_{ij}(0),\dot\gamma_{jk}(0))>0$ because 
$\pi_j(p)=\Span_\R(\dot\gamma_{ij}(0),\dot\gamma_{jk}(0))$ is a symplectic 
plane. Applying a symplectic linear change of coordinates and a rescaling of the 
parameterizations, we can suppose that 
$$
\dot\gamma_{12}(0)=
\left(
\begin{array}{l}
1 \\0 \\1\\0
\end{array}\right)\, , \;\;
\dot\gamma_{23}(0)=
\left(
\begin{array}{l}
0 \\1 \\1\\0
\end{array}\right)\, , \;\;
\dot\gamma_{31}(0)=
\left(
\begin{array}{l}
-1 \\-1 \\1\\0
\end{array}\right)\,\subset \C^2\approx \R^4.
$$ 
Looking sufficiently close to $p$, there is a diffeomorphism $\Phi$ $\cc^1$-close to 
the identity taking $S$ to the cone over the triangle :
$$
\Sigma:=\left\{ \left(\begin{array}{l}
z\phi_1(\theta)\\z \\0
\end{array}\right)\, ,\;0<h<1\right\}
$$
where $\phi_1(\theta)$ is the counterclockwise parameterization of the triangle 
$T_1$ spanned by the points $(1,0)$, $(0,1)$, $(-1,-1)$ of $\R^2$. 
Consider then the polar parameterizations 
$\phi_\rho$ of a family of convex curves $T_\rho$ with $T_\rho=\{x^2+y^2=\rho^4\}$ 
for $\rho\ll 1$, $T_\rho=\rho\cdot T_1$ for 
$\rho\simeq 1$ and the angles of the singularities of $T_\rho$ vary smoothly with 
$\rho$ (see figure \ref{ouf}). Then the surface
$$
\wdt \sigma:=\{(\phi_h(\theta),h,0)\}
$$
is smooth near the origin, has tangent planes far from being Lagrangian and coincide
with $\Sigma$ in a neighbourhood of $h=\{y_2=1\}$. Cutting $\Sigma$ from $S$ and 
replacing it by $\Phi^{-1}(\wdt\Sigma)$, we get a symplectic smoothening of $S$ at $p$.
We were also sufficiently cautious in the extrapolation from $T_1$ to $T_\rho$, 
$\rho\simeq 0$ to ensure that the 
remaining singularities have the form (\ref{nicecyl}).\hfill $\square$

{  
\footnotesize
\bibliography{bib3.bib}
\bibliographystyle{abbrv}
}
\vspace{5cm}

\noindent Emmanuel Opshtein,\\
\noindent École normale supérieure de Lyon\\
\noindent Unité de Mathématiques Pures et Appliquées\\
\noindent UMR CNRS 5669\\
\noindent 46, allée d'Italie\\
\noindent 69364 LYON Cedex 07\\
\noindent FRANCE  \vspace{1cm}\\
\noindent E-mail : emmanuel.opshtein@umpa.ens-lyon.fr
\end{document}

%% file: karshon.pstex_t
\begin{picture}(0,0)%
\includegraphics{karshon.pstex}%
\end{picture}%
\setlength{\unitlength}{2693sp}%
\begingroup\makeatletter\ifx\SetFigFont\undefined%
\gdef\SetFigFont#1#2#3#4#5{%
  \reset@font\fontsize{#1}{#2pt}%
  \fontfamily{#3}\fontseries{#4}\fontshape{#5}%
  \selectfont}%
\fi\endgroup%
\begin{picture}(6286,2253)(1201,-2605)
\put(1306,-511){\makebox(0,0)[lb]{\smash{{\SetFigFont{9}{10.8}{\rmdefault}{\mddefault}{\updefault}{\color[rgb]{0,0,0}$A$}%
}}}}
\put(3106,-2491){\makebox(0,0)[lb]{\smash{{\SetFigFont{9}{10.8}{\rmdefault}{\mddefault}{\updefault}{\color[rgb]{0,0,0}$C$}%
}}}}
\put(1216,-1501){\makebox(0,0)[lb]{\smash{{\SetFigFont{9}{10.8}{\rmdefault}{\mddefault}{\updefault}{\color[rgb]{0,0,0}$I$}%
}}}}
\put(1261,-2536){\makebox(0,0)[lb]{\smash{{\SetFigFont{9}{10.8}{\rmdefault}{\mddefault}{\updefault}{\color[rgb]{0,0,0}$B$}%
}}}}
\put(5581,-1276){\makebox(0,0)[lb]{\smash{{\SetFigFont{9}{10.8}{\rmdefault}{\mddefault}{\updefault}{\color[rgb]{0,0,0}$I$}%
}}}}
\put(5581,-2536){\makebox(0,0)[lb]{\smash{{\SetFigFont{9}{10.8}{\rmdefault}{\mddefault}{\updefault}{\color[rgb]{0,0,0}$B$}%
}}}}
\put(7426,-2491){\makebox(0,0)[lb]{\smash{{\SetFigFont{9}{10.8}{\rmdefault}{\mddefault}{\updefault}{\color[rgb]{0,0,0}$C$}%
}}}}
\put(1261,-2536){\makebox(0,0)[lb]{\smash{{\SetFigFont{9}{10.8}{\rmdefault}{\mddefault}{\updefault}{\color[rgb]{0,0,0}$B$}%
}}}}
\put(5626,-511){\makebox(0,0)[lb]{\smash{{\SetFigFont{9}{10.8}{\rmdefault}{\mddefault}{\updefault}{\color[rgb]{0,0,0}$A$}%
}}}}
\end{picture}%

%% file: pack5.pstex_t
\begin{picture}(0,0)%
\includegraphics{pack5.pstex}%
\end{picture}%
\setlength{\unitlength}{2652sp}%
\begingroup\makeatletter\ifx\SetFigFont\undefined%
\gdef\SetFigFont#1#2#3#4#5{%
  \reset@font\fontsize{#1}{#2pt}%
  \fontfamily{#3}\fontseries{#4}\fontshape{#5}%
  \selectfont}%
\fi\endgroup%
\begin{picture}(10317,3262)(719,-2693)
\put(2216,-811){\makebox(0,0)[lb]{\smash{\SetFigFont{9}{10.8}{\rmdefault}{\mddefault}{\updefault}{\color[rgb]{0,0,0}$p_k$}%
}}}
\put(5896,389){\makebox(0,0)[lb]{\smash{\SetFigFont{9}{10.8}{\rmdefault}{\mddefault}{\updefault}{\color[rgb]{0,0,0}$\Psi$}%
}}}
\put(1406,-1391){\makebox(0,0)[lb]{\smash{\SetFigFont{9}{10.8}{\rmdefault}{\mddefault}{\updefault}{\color[rgb]{0,0,0}$p_2$}%
}}}
\put(3826,-2626){\makebox(0,0)[lb]{\smash{\SetFigFont{9}{10.8}{\rmdefault}{\mddefault}{\updefault}{\color[rgb]{0,0,0}$(\D\times \D,\om')$}%
}}}
\put(2429,-661){\makebox(0,0)[lb]{\smash{\SetFigFont{9}{10.8}{\rmdefault}{\mddefault}{\updefault}{\color[rgb]{0,0,0}$(\phi(z_1),z_2)$}%
}}}
\put(2429,-1161){\makebox(0,0)[lb]{\smash{\SetFigFont{9}{10.8}{\rmdefault}{\mddefault}{\updefault}{\color[rgb]{0,0,0}$\phi^*\sigma_{st}=\sigma_{st}$}%
}}}
\put(5573,-661){\makebox(0,0)[lb]{\smash{\SetFigFont{9}{10.8}{\rmdefault}{\mddefault}{\updefault}{\color[rgb]{0,0,0}$\Phi$}%
}}}
\put(7287,-375){\makebox(0,0)[lb]{\smash{\SetFigFont{9}{10.8}{\rmdefault}{\mddefault}{\updefault}{\color[rgb]{0,0,0}$B_R$}%
}}}
\put(786,-2589){\makebox(0,0)[lb]{\smash{\SetFigFont{9}{10.8}{\rmdefault}{\mddefault}{\updefault}{\color[rgb]{0,0,0}$(D\times \D,\om)$}%
}}}
\put(6796,-2626){\makebox(0,0)[lb]{\smash{\SetFigFont{9}{10.8}{\rmdefault}{\mddefault}{\updefault}{\color[rgb]{0,0,0}$(\ce_{R,1/\sqrt 2},\om_{st})$}%
}}}
\put(9001,-2626){\makebox(0,0)[lb]{\smash{\SetFigFont{9}{10.8}{\rmdefault}{\mddefault}{\updefault}{\color[rgb]{0,0,0}$\Psi^{-1}(B_R)$}%
}}}
\put(9446,-1916){\makebox(0,0)[lb]{\smash{\SetFigFont{9}{10.8}{\rmdefault}{\mddefault}{\updefault}{\color[rgb]{0,0,0}$p_2$}%
}}}
\put(11036,-931){\makebox(0,0)[lb]{\smash{\SetFigFont{9}{10.8}{\rmdefault}{\mddefault}{\updefault}{\color[rgb]{0,0,0}$p_k$}%
}}}
\put(8324,-966){\makebox(0,0)[lb]{\smash{\SetFigFont{9}{10.8}{\rmdefault}{\mddefault}{\updefault}{\color[rgb]{0,0,0}$p_1$}%
}}}
\put(719,-891){\makebox(0,0)[lb]{\smash{\SetFigFont{9}{10.8}{\rmdefault}{\mddefault}{\updefault}{\color[rgb]{0,0,0}$p_1$}%
}}}
\end{picture}

%% file: fun3.pstex_t
\begin{picture}(0,0)%
\includegraphics{fun3.pstex}%
\end{picture}%
\setlength{\unitlength}{1367sp}%
\begingroup\makeatletter\ifx\SetFigFont\undefined%
\gdef\SetFigFont#1#2#3#4#5{%
  \reset@font\fontsize{#1}{#2pt}%
  \fontfamily{#3}\fontseries{#4}\fontshape{#5}%
  \selectfont}%
\fi\endgroup%
\begin{picture}(13299,4290)(-1586,-5041)
\put(6661,-4666){\makebox(0,0)[lb]{\smash{\SetFigFont{5}{6.0}{\rmdefault}{\mddefault}{\updefault}{\color[rgb]{0,0,0}$\ell_r$}%
}}}
\put(2394,-2318){\makebox(0,0)[lb]{\smash{\SetFigFont{5}{6.0}{\rmdefault}{\mddefault}{\updefault}{\color[rgb]{0,0,0}$f$}%
}}}
\put(-854,-3706){\makebox(0,0)[lb]{\smash{\SetFigFont{5}{6.0}{\rmdefault}{\mddefault}{\updefault}{\color[rgb]{0,0,0}$\pi_1(p)$}%
}}}
\put(8551,-1696){\makebox(0,0)[lb]{\smash{\SetFigFont{5}{6.0}{\rmdefault}{\mddefault}{\updefault}{\color[rgb]{0,0,0}$\ell_l$}%
}}}
\put(4846,-1013){\makebox(0,0)[lb]{\smash{\SetFigFont{5}{6.0}{\rmdefault}{\mddefault}{\updefault}{\color[rgb]{0,0,0}$2\pi$}%
}}}
\put(11251,-2986){\makebox(0,0)[lb]{\smash{\SetFigFont{5}{6.0}{\rmdefault}{\mddefault}{\updefault}{\color[rgb]{0,0,0}$T_p\ell_l=T_p\ell_r$}%
}}}
\put(1126,-1651){\makebox(0,0)[lb]{\smash{\SetFigFont{5}{6.0}{\rmdefault}{\mddefault}{\updefault}{\color[rgb]{0,0,0}$\pi_2(p)$}%
}}}
\put(3541,-4186){\makebox(0,0)[lb]{\smash{\SetFigFont{6}{7.2}{\rmdefault}{\mddefault}{\updefault}{\color[rgb]{0,0,0}$Q_\eps$}%
}}}
\put(5266,-2851){\makebox(0,0)[lb]{\smash{\SetFigFont{5}{6.0}{\rmdefault}{\mddefault}{\updefault}{\color[rgb]{0,0,0}$p$}%
}}}
\end{picture}

%% file: ouf.pstex_t
\begin{picture}(0,0)%
\includegraphics{ouf.pstex}%
\end{picture}%
\setlength{\unitlength}{912sp}%
\begingroup\makeatletter\ifx\SetFigFont\undefined%
\gdef\SetFigFont#1#2#3#4#5{%
  \reset@font\fontsize{#1}{#2pt}%
  \fontfamily{#3}\fontseries{#4}\fontshape{#5}%
  \selectfont}%
\fi\endgroup%
\begin{picture}(20744,7599)(-12171,-7732)
\put(6346,-3571){\makebox(0,0)[lb]{\smash{\SetFigFont{5}{6.0}{\rmdefault}{\mddefault}{\updefault}{\color[rgb]{0,0,0}$\theta$}%
}}}
\put(5086,-1096){\makebox(0,0)[lb]{\smash{\SetFigFont{5}{6.0}{\rmdefault}{\mddefault}{\updefault}{\color[rgb]{0,0,0}$T_1$}%
}}}
\put(-8574,-5371){\makebox(0,0)[lb]{\smash{\SetFigFont{5}{6.0}{\rmdefault}{\mddefault}{\updefault}{\color[rgb]{0,0,0}$\pi_1(p)$}%
}}}
\put(-6704,-3886){\makebox(0,0)[lb]{\smash{\SetFigFont{5}{6.0}{\rmdefault}{\mddefault}{\updefault}{\color[rgb]{0,0,0}$\pi_3(p)$}%
}}}
\put(-8084,-7651){\makebox(0,0)[lb]{\smash{\SetFigFont{5}{6.0}{\rmdefault}{\mddefault}{\updefault}{\color[rgb]{0,0,0}$p$}%
}}}
\put(-11924,-736){\makebox(0,0)[lb]{\smash{\SetFigFont{5}{6.0}{\rmdefault}{\mddefault}{\updefault}{\color[rgb]{0,0,0}$\ell_{12}$}%
}}}
\put(-6149,-496){\makebox(0,0)[lb]{\smash{\SetFigFont{5}{6.0}{\rmdefault}{\mddefault}{\updefault}{\color[rgb]{0,0,0}$\ell_{23}$}%
}}}
\put(-3374,-706){\makebox(0,0)[lb]{\smash{\SetFigFont{5}{6.0}{\rmdefault}{\mddefault}{\updefault}{\color[rgb]{0,0,0}$\ell_{13}$}%
}}}
\put(-9565,-3732){\makebox(0,0)[lb]{\smash{\SetFigFont{5}{6.0}{\rmdefault}{\mddefault}{\updefault}{\color[rgb]{0,0,0}$\pi_2(p)$}%
}}}
\end{picture}